\documentclass[]{amsproc}

\usepackage{verbatim}
\usepackage{amssymb}
\usepackage{amsbsy}
\usepackage{amscd}
\usepackage{amsmath}
\usepackage{amsthm}
\usepackage[mathscr]{eucal}

\newtheorem{theorem}{Theorem}
\newtheorem{defn}{Definition}\numberwithin{defn}{section}
\newtheorem{conjecture}{Conjecture}
\newtheorem{cor}{Corollary}
\newtheorem{prop}{Proposition}

\title[On a conjecture by N. Takahashi]{On a conjecture by N. Takahashi on log mirror symmetry for the projective plane}

\author[M. van Garrel]{Michel van Garrel}

\address{KIAS, 85 Hoegiro Dongdaemun-gu, Seoul 130-722, Republic of Korea}

\email{vangarrel@kias.re.kr}

\subjclass[2010]{05A15, 14J33, 14J45, 14N35}

\keywords{}

\newcommand{\ptwo}{\mathbb{P}^2}
\newcommand{\pone}{\mathbb{P}^1}

\newcommand{\iz}{\mathbb Z}
\newcommand{\iq}{\mathbb Q}

\DeclareMathOperator{\hhh}{H}
\DeclareMathOperator{\nnn}{N}
\DeclareMathOperator{\ttt}{T}
\DeclareMathOperator{\euler}{e}
\DeclareMathOperator{\vdim}{vdim}

\begin{document}

\begin{abstract}

In this survey, the relationship between local and relative BPS state counts is explored in light of a conjecture on log mirror symmetry by N. Takahashi. This is based on the paper \cite{GWZ} by the author and the integrality result of \cite{GWZ} is reported.

\end{abstract}

\maketitle

\section{Introduction}

Log mirror symmetry for the pair of projective plane and elliptic curve on it was explored by N. Takahashi in \cite{Ta}. The $A$-model invariants considered are counts of rational curves in the plane meeting the elliptic curve in a prescribed point of prescribed order. A relationship of these invariants with local BPS state counts is conjectured, the latter being the $A$-model invariants of local $\ptwo$. In the author's paper \cite{GWZ}, a virtual version of Takahashi's curve counts, derived from relative Gromov-Witten invariants, are considered. Namely, the relative BPS state counts introduced by Gross-Pandharipande-Siebert in \cite{GPS}. The main result of \cite{GWZ} is the integrality of these invariants for toric del Pezzo surfaces. The key result towards proving that integrality claim is a linear relationship between the local and relative BPS state counts. Incidentally, this provides a virtual reformulation of Takahashi's conjecture, which we describe in the present survey. We start in section \ref{prelim} by defining local BPS state counts, putting extra care into justifying the definition. We then introduce Takahashi's conjecture in section \ref{taka}. In \cite{Ta}, the conjectural $B$-model to log mirror symmetry is described. Namely, the mirror family and the associated periods are defined. We do not review them as our correspondence of BPS numbers concerns the $A$-model. Finally, in section \ref{relbps} we define relative BPS state counts and explain how they satisfy an analogue to Takahashi's conjecture.

\section{Preliminaries}
\label{prelim}

Let $S$ be a smooth del Pezzo surface. In definition \ref{localbps} below, local BPS state counts of $S$ are introduced. These are the $A$-model invariants of local mirror symmetry, which was developed in \cite{CKYZ}. See also \cite{KM} for a description in terms of Yukawa couplings.

Denote moreover by $D$ a smooth effective anticanonical divisor on $S$, by $K_S$ the total space of the canonical bundle $\mathcal{O}_S(-D)$, and let $\beta\in\hhh_{2}(S,\iz)$. Denote by $\overline{M}_{0,0}(S,\beta)$, resp. by $\overline{M}_{0,1}(S,\beta)$, the moduli stack of genus 0 stable maps
$$
f:C\to S
$$
with no, resp. one, marked point and such that $f_{*}([C])=\beta$. Denote by
$$
\pi:\overline{M}_{0,1}(S,\beta)\to\overline{M}_{0,0}(S,\beta)
$$
the forgetful morphism and by
$$
ev:\overline{M}_{0,1}(S,\beta)\to S
$$
the evaluation map. This determines the bundle $R^{1}\pi_{*}ev^{*}K_{S}$ whose fiber over a stable map $f:C\to S$ is $\hhh^{1}(C,f^{*}K_{S})$. Consider moreover the virtual fundamental class $[\overline{M}_{0,0}(S,\beta)]^{vir}$, which intersection-theoretically is the correct class to be looking at, and denote by $\euler$ the Euler class.

The idea of obtaining \emph{local} Gromov-Witten invariants of $S$ comes from assuming that $S$ was embedded in a compact Calabi-Yau threefold $X$. That is, $X$ is a 3-dimensional K\"ahler manifold with trivial canonical bundle. If $S\subset X$, then the local Gromov-Witten invariants of $S$ should be the contribution of stable maps mapping to $S$ to the Gromov-Witten invariants of $X$. Through excess intersection calculations, the authors of \cite{CKYZ} arrive at the following definition.

\begin{defn}
The \emph{genus 0 degree $\beta$ local Gromov-Witten invariant of $S$} is
$$
I_{K_{S}}(\beta):=\int_{[\overline{M}_{0,0}(S,\beta)]^{vir}}\euler\left(R^{1}\pi_{*}ev^{*}K_{S}\right)\in\mathbb{Q}.
$$
\end{defn}

The degree map on a proper Deligne-Mumford stack is $\iq$-valued and thus in general, the invariants $I_{K_S}(\beta)$ are $\iq$-valued. This makes a direct enumerative interpretation impossible. Local BPS state counts aim at remedying this by extracting integer-valued invariants from the rational number $I_{K_S}(\beta)$. These new invariants would under some idealised (usually not satisfied) conditions be counts of rational curves.

As above, assume that $S$ is embedded in a compact Calabi-Yau threefold $X$. Then $\nnn_{S/X} = K_S$, so that $K_S$ is the local geometry of $X$ near $S$. Moreover, the $I_{K_S}(\beta)$ can be regarded as counts of maps into $K_S$ since genus 0 stable maps $C \to K_S$ factor through $S$. Indeed, for simplicity assume that $C=\pone$. Then giving a map $\phi : \pone \to K_S$ is the same as a map $f : \pone \to S$ and a section of $f^{*}\left( \mathcal{O}(-D) \right)$. Since $S$ is del Pezzo however, the latter is negative on $\pone$ and thus has no sections. It follows that the image of $\pone$ in $K_S$ is contained in the 0-section. Therefore, the moduli space of stable genus 0 maps into $K_S$ is identical to the moduli space of stable genus 0 maps into $S$. However, their deformation theories differ yielding different virtual fundamental classes. For one, the virtual fundamental class of maps into $K_S$ is in degree 0, whereas it is in degree non-zero for $S$.

Via $S \subseteq X$, we model stable maps to $S$ from the situation for $X$, which we describe now. Let $\phi : \pone \to X$ be a holomorphic immersion. That is, we are looking at a rational curve $\phi(\pone)$ in $X$ (with possibly nodes), and $\phi$ is the normalization map. Then the normal sheaf of $\phi$ is
$$
\nnn_\phi =  \dfrac{\phi^*\ttt_X}{\ttt_{\pone}} \cong \mathcal{O}_{\pone}(a)\oplus\mathcal{O}_{\pone}(b),
$$
where by the adjunction formula $a+b=-2$. Moreover, $a=b=-1$ if and only if $\nnn_\phi$ has no sections if and only if $\phi(\pone)$ can not be infinitesimally deformed. Hence, $\phi(\pone)$ is infinitesimally rigid if and only if $a=b=-1$.

We come to the Aspinwall-Morrison formula, which was proved in \cite{Ma}, see also \cite{Vo}: Assume that $\phi(\pone)$ is infinitesimally rigid of degree $\beta$, so that $\nnn_\phi \cong \mathcal{O}_{\pone}(-1)\oplus\mathcal{O}_{\pone}(-1)$. Consider maps $\pone \to X$ with image $\phi(\pone)$. Such a map factors through the normalization of $\phi(\pone)$, i.e. we are looking at multiple covers of $\pone$. For the purpose of calculating Gromov-Witten invariants, we are considering multiple covers of $\pone$ into $\nnn_\phi$. The Aspinwall-Morrison formula calculates the contributions of such maps to the Gromov-Witten invariants of $X$. It states that the degree $k$ covers
$$
\pone \xrightarrow{k:1} \pone,
$$
and their degenerations, (after composing with $\phi$) contribute $\frac{1}{k^3}$ to the Gromov-Witten invariant $I_X(k\beta)$.

We come to the three idealised assumptions that underlie the definition of BPS numbers:
\begin{itemize}
\item Each rational curve in $X$ is infinitesimally rigid.
\item Any two irreducible rational curves do not intersect.
\item For each curve class $\beta$, the number of infinitesimally rigid rational curves of class $\beta$ is finite.
\end{itemize}
Supposing that these three idealised assumptions are satisfied, denote by $n'_\beta$ the number of rational curves in $X$ of degree $\beta$. Let $\gamma$ be a curve class. Because of our assumptions, the only way to obtain a stable map of degree $\gamma$ is as follows: Decompose $\gamma = k \cdot \beta$, for $k$ a positive integer and $\beta$ a curve class. Then consider the maps consisting of a degree $k$ cover of $\pone$ followed by a holomorphic immersion to a rational curve of degree $\beta$ in $X$. By the Aspinwall-Morrison formula, all possible such decompositions contribute to $I_X(\gamma)$ as follows:
$$
I_X(\gamma) = \sum_{\gamma = k\cdot\beta}\frac{1}{k^3} \; n'_\beta.
$$
In particular, this stipulates that the rational numbers of Gromov-Witten invariants come from multiple covers. 

The preceding heuristic argument justifies the definition below of local BPS state counts. Note that in absence of the prior mentioned ideal conditions, the meaning of the $n_\beta$ is no longer clear. Moreover, the previous discussion assumes a compact Calabi-Yau threefold, whereas we are dealing with the open Calabi-Yau threefold $K_S$. Nonetheless, the $n_\beta$ below are well-defined. They are the $A$-model invariants of local mirror symmetry for $S$, which is why they are crucial.

\begin{defn}\label{localbps} (See \cite{GVI, GVII, PH, P3}.)
Assume (for simplicity of exposition) that $\beta$ is primitive and let $d\geq 1$. Then the \emph{local BPS state counts}  $n_{d\beta}\in\iq$ are defined via the equality of power series:
$$
\sum_{l=1}^{\infty}I_{K_{S}}(l\beta)\, q^{l}=\sum_{d=1}^{\infty}n_{d\beta}\sum_{k=1}^{\infty}\frac{1}{k^{3}}\, q^{dk}.
$$
\end{defn}

Note that vice versa, one can pass from the $n_\beta$ to the $I_{K_S}(\beta)$ by a M\"obius transform. The following theorem was originally a conjecture attributed to Gopakumar-Vafa and stated in \cite{BPB}. It was proved in \cite{Pe} in the case of toric del Pezzo surfaces. These are the del Pezzo surfaces of degree $\geq 6$. More generally for toric Calabi-Yau threefolds, a proof was given in \cite{Ko}.

\begin{theorem}\label{localint}(Peng in \cite{Pe}.) If $S$ is a toric del Pezzo surface, then
$$
n_{d\beta}\in\iz
$$
for all $d\geq 1$.
\end{theorem}

We finish our preliminaries with the definition of the geometry that will be relevant in the next sections.

\begin{defn}\label{logCY} (See \cite{GPS}.)
Denote by $S$ a smooth surface, by $D\subseteq S$ a smooth divisor and let $0\neq\gamma\in\hhh_2(S,\iz)$. The pair $(S,D)$ is said to be \emph{log Calabi-Yau with respect to $\gamma$} if
$$
D\cdot\gamma = c_1(S)\cdot\gamma.
$$
Provided the above equation holds for all $\gamma$, $(S,D)$ is said to be \emph{log Calabi-Yau}.
\end{defn}

A class of examples of log Calabi-Yau surface pairs are the $(S,D)$ for any del Pezzo surface $S$ and any smooth anti-canonical curve $D$ on it.

\section{$A$-model log mirror symmetry for $\ptwo$}
\label{taka}

Log mirror symmetry for the projective plane was explored by N. Takahashi in \cite{Ta}. Denote by $E\subseteq\ptwo$ a smooth cubic, i.e. an elliptic curve, so that $(\ptwo,E)$ is a log Calabi-Yau surface pair. According to the conjecture of \cite{Ta}, the log mirror symmetry $A$-model concerns counts of plane curves satisfying \emph{condition (AL)}, which stands for \emph{affine line}:

\begin{defn}
Let $C\subseteq\ptwo$ be a curve. We say that $C$ is an \emph{(AL)-curve} if $C$ is irreducible, reduced and if the normalization of $C\setminus E$ is isomorphic to the affine line $\mathbb{A}^1$.
\end{defn}

Let $C$ be an (AL)-curve of degree $d$. Then by definition, $C$ is of geometric genus 0 and meets $E$ in exactly one point. If $P$ is that point, then the equation defining $C$ yields the linear equivalence $3d\cdot P \sim 3d \cdot P_0$, where $P_0$ is any one of the nine inflexion points of $E$. Recall that the identification of closed points of $E$ with its divisor class group yields, upon choosing a zero element, a group law on the set of closed points of $E$. Hence by choosing as zero element $0\in E$ a flex point, it follows that $P$ is a $3d$-torsion point for the resulting group law. Takahashi's proposal for $A$-model invariants is to count the \emph{(AL)}-curves that meet $E$ in a chosen point of maximal order:

\begin{defn}
Choose $0\in E$ to be a flex point and let $P$ be a point of order $3d$ for the group law induced by the above choice of $0$. For $d\geq 1$, define $m_d$ to be the number of \emph{(AL)}-curves of degree $d$ that meet $E$ in $P$.
\end{defn}

In \cite{Ta}, it is conjectured that the $B$-model to $(\ptwo, E)$ is the same as the $B$-model to local $\ptwo$. We refer to \cite{Ta} for the precise statements. Instead, we focus on the conjectured relationship of the $A$-model invariants $m_d$ to the $A$-model invariants of local $\ptwo$. An equivalent way of stating Takahashi's conjecture then is as follows. Denote by $n_d$ the local BPS state counts of $\ptwo$, i.e. we suppress the hyperplane class from the notation.

\begin{conjecture}\label{conj} (N. Takahashi in \cite{Ta}.) For all $d\geq 1$,
$$
3d \, m_d = (-1)^{d+1}  n_d.
$$
\end{conjecture}

In light of theorem \ref{thm}, a more natural statement of the above conjecture is $9d^2 \, m_d = (-1)^{d+1} \, 3d \, n_d$. Note that $9d^2$ is the number of $3d$-torsion points on $E$.

We remarked in section \ref{prelim} that an enumerative interpretation of the local BPS state counts remains elusive. The above conjecture would however provide exactly that. According to it, the $n_d$ are counts of rational curves in the projective plane meeting an elliptic curve in a prescribed point of order $3d$.

\section{Relative BPS state counts of del Pezzo surfaces}
\label{relbps}

The pair $(S,D)$ can be thought of as being log Calabi-Yau because the canonical bundle of $S$ is trivial away from the divisor $D$. The definition of relative BPS state counts mirrors the definition of local BPS state counts of section \ref{prelim}. While we were considering stable maps then, we now consider relative stable maps, i.e. we additionally prescribe how the maps meet $D$. Let $\beta\in\hhh_{2}(S,\iz)$ be the class of a curve and set $w=D\cdot\beta$. Then a generic stable map of degree $\beta$ meets $D$ in $w$ points. Moreover, the moduli space of stable maps of degree $\beta$ is of virtual dimension
$$
\int_\beta c_1(\ttt_S) - 1 =  w-1.
$$ 
Each time two points of intersections are identified, the virtual dimension drops by one. Proceeding accordingly with all but one of the intersection points therefore cuts it by $w-1$ and hence, in the minimal case of one intersection point only, yields a virtual dimension of 0.

The virtual counts of such maps are given by relative Gromov-Witten invariants, which were introduced to the algebraic setting by Li in \cite{Li}. Denote by $\overline{M}(S/D,w)$ the moduli space of relative genus 0 stable maps of degree $\beta$ meeting $D$ in one point of maximal tangency $w$. Since
$$
\vdim\overline{M}(S/D,w)=0,
$$
the following definition makes sense.
\begin{defn}
The \emph{genus 0 degree $\beta$ relative Gromov-Witten invariant of maximal tangency} is the degree of the corresponding virtual fundamental class:
$$
N_{S}[w]:=\int_{[\overline{M}(S/D,w)]^{vir}}1\in\mathbb{Q}.
$$
\end{defn}
Note that for notational simplicity, $\beta$ is hidden in the notation. Analogously to the local case, BPS numbers are defined via considering multiple cover contributions over rigid elements. Denote by $\iota:P\to S$ a rigid element of $\overline{M}(S/D,w)$ and for $k\geq1$, let $M_{P}[k]$ be the contribution of $k$-fold multiple covers of $P$ to $N_{S}[kw]$ (see \cite{GPS} for precise definitions).
\begin{prop}
(Proposition $6.1$ in \cite{GPS}.)
\[
M_{P}[k]=\frac{1}{k^{2}}\binom{k(w-1)-1}{k-1}.
\]
\end{prop}
The next definition exactly mirrors the heuristic discussion of section \ref{prelim} with corresponding idealised assumptions for the relative case of one point of intersection.
\begin{defn}
\label{def-rel-bps}(Paragraph 6.3 in \cite{GPS}.)
Let $d\geq1$. Then the \emph{relative BPS state counts} $n_{S}[dw]\in\iq$ are
defined by means of the equality
\begin{equation}
\sum_{l=1}^{\infty}N_{S}[lw]\, q^{l}=\sum_{d=1}^{\infty}n_{S}[dw]\sum_{k=1}^{\infty}\frac{1}{k^{2}}\binom{k(dw-1)-1}{k-1}\, q^{dk}.\label{eq:defn_rel_bps_gen_fn}
\end{equation}
\end{defn}

As in the local case, a M\"obius transform recovers the $N_S[dw]$ from the $n_S[dw]$. Analogously to the local invariants, the relative BPS state counts are conjectured to be integers:

\begin{conjecture}
\label{conjecture_GPS}(Conjecture $6.2$ in \cite{GPS}.)
For $\beta\in\hhh_{2}(S,\iz)$ a primitive curve class, set
$w=\beta\cdot D$. Then, for all $d\geq1$,
\[
n_{S}[dw]\in\iz.
\]
\end{conjecture}

In the idealised situation where every rational curve meeting $D$ in one point was rigid, and where there are only finitely many such curves in each degree, the relative BPS state counts would correspond to the (AL)-curve counts, except that the point of intersection with $D$ is not prescribed. In that sense, the $n_S[dw]$ are a double generalisation of the $m_d$: Firstly, any point of any order (for the group law) is allowed as intersection point. Secondly, we consider virtual invariants, i.e. we take the degree of the virtual fundamental class of the relevant moduli stack.

Before proceeding with the linear restatement of conjecture \ref{conj}, we introduce some notation. For  $n\in\iz$, denote by $\omega(n)$ the number of primes (not
counting multiplicities) in the prime factorization of $n$. Furthermore, let
$$
I(n):=\left\{ k\in\iz\,:\, k|n\text{ and }n/k\text{ is square-free}\right\}.
$$
We now introduce the following infinite-dimensional matrix $C$. If $t|s$, set
$$
C_{st}:=\frac{(-1)^{sw}}{(s/t)^{2}}\sum_{k\in I(s/t)}(-1)^{\omega\left(s/kt\right)}(-1)^{ktw}\binom{k(tw-1)-1}{k-1}.
$$
If $t\nmid s$, set $C_{st}=0$. Consequently, each row of $C$ has only a finite number of non-zero entries. Therefore, we can apply to $C$ an infinite-dimensional vector without encountering convergence issues.

We are ready to state a virtual version of conjecture \cite{Ta} by N. Takahashi. Our result applies to all del Pezzo surfaces.

\begin{theorem}\label{thm} (See lemma 12 in \cite{GWZ}.)
$$
C\cdot\left[n_{S}[dw]\right]_{d\geq 1}=\left[(-1)^{dw+1}\, dw\, n_{d\beta}\right]_{d\geq 1}.
$$
\end{theorem}

In \cite{GWZ}, it is proven that the entries of $C$ are integers. Since $\det(C)=1$, the entries of $C^{-1}$ are integers as well. Hence, the following integrality result follows from theorem \ref{localint}.

\begin{cor}\label{cor} (See Corollary 10 in \cite{GWZ}.)
For toric del Pezzo surfaces, the relative BPS state counts are integers.
\end{cor}

%%%%%%%%%%%%%%%%%%%%%%%%%%%%%%%%%
% References
%%%%%%%%%%%%%%%%%%%%%%%%%%%%%%%%%

\end{document}